\documentclass[12pt,reqno]{amsart}
\usepackage{amsmath,amsopn,amssymb,amsthm,multicol}
\usepackage{hyperref}
\usepackage{graphicx}

\newcommand{\fr}{\mathfrak}

 \newtheorem{lemma} {Lemma} [section]
\newtheorem{theorem}[lemma]{Theorem}

\newtheorem{corol}[lemma] {Corollary} 
\newtheorem{example}[lemma] {Example}

\begin{document}

\title{Motion of charged particle in a class of homogeneous spaces}
\author{Andreas Arvanitoyeorgos  and Nikolaos Panagiotis Souris}\thanks{The second author was partially supported by the EPSRC grant EP/R008205/1.  Part of this work was completed during his postdoctoral position at the University of Reading, UK. An older
version of this work was supported by Grant E.037 from the Research Committee of the University
of Patras (Programme K. Karatheodori).} 
\address{University of Patras, Department of Mathematics, GR-26500 Patras, Greece}
\email{arvanito@math.upatras.gr}

\address{University of Reading, Department of Mathematics and Statistics, Whiteknights, PO Box 220, Reading RG6 6AX, UK}
\email{ntgeva@hotmail.com}
\medskip

\begin{abstract} 
We study the motion of charged particle under a natural choice of electromagnetic field in a general class of compact homogeneous spaces.  As a special case we describe the motion in homogeneous Riemannian spaces $(G/H,g)$, where $g$ is any deformation of a normal metric along the fibers of a homogeneous fibration $K/H\rightarrow G/H\rightarrow G/K$.     


\end{abstract}

\maketitle
\medskip

\section{Introduction and statement of results}

Let $(M=G/H, g)$ be a homogeneous Riemannian manifold, where $g$ is a $G$-invariant metric with corresponding Levi-Civita connection $\nabla$.  
The aim of the present paper is the study of the charged particle differential equation
\begin{equation}\label{cp}
\nabla_{\dot{x}}{\dot{x}}=kI(\dot{x}),
\end{equation} 
where $k$ is some real constant and $I$ a (1,1) tensor in $M$ that corresponds to an electromagnetic field.  The curve $x(t)$, describing the motion of a charged particle, is also called a \emph{magnetic curve}.
Equation (\ref{cp}) appears in a more general context in general relativity as follows (\cite{MTW}).

Let $(M,g)$ be a Riemannian manifold, $F$ a closed $2$-form, and $X$ a
vector field on $M$.  We denote by $\iota _X:\Lambda ^p(M)\to\Lambda
^{p-1}(M)$ the interior product operator induced by $X$, and by
$\mathcal{L}: TM\to T^*M$ the Legendre transformation defined by
$u\mapsto\mathcal{L}(u)$, $\mathcal{L}(u)(v)=g(u,v)$ ($v\in TM$). A
curve $x(t)$ in $M$ is called a {\it motion of a charged particle
under electromagnetic field $F$} if it satisfies the differential
equation
$$
\nabla _{\dot{x}}\dot{x}=-\mathcal{L}^{-1}(\iota _{\dot{x}}F),
$$
where $\nabla$ is the Levi-Civita connection of $M$.  When $F=0$ then
$x(t)$ is a geodesic in $M$. 
In particular, if 
 $M$ is a K\"ahler manifold with complex structure $J$ there is a
natural choice of an electromagnetic field $F$, namely a scalar
multiple of the K\"ahler form $\omega$, defined by $\omega
(X,Y)=g(X, JY)$.  Since $-\mathcal{L}^{-1}(\iota _X\omega)=JX$, a
curve $x(t)$ is a motion of a charged particle under
electromagnetic field $\kappa\omega$ if and only if
$
\nabla _{\dot{x}}\dot{x}=kJ(\dot{x})$.
We also refer to \cite{KST} and \cite{Li} for  other relevant applications in physics.

Magnetic curves are being extensively studied in various classes of spaces (see for example \cite{Bo-Jo}, \cite{CFG}, \cite{Er-In}, \cite{Ik3}, \cite{MN}, and references therein).  
A class of homogeneous spaces considered by Ikawa in \cite{Ik3} were generalized flag manifolds 
(or K\"ahler C-spaces) with two isotropy summands.
These spaces were classified by the first author and I. Chrysikos in \cite{Arv-Chr}, hence obtaining a concrete class of homogeneous spaces where equation (\ref{cp}) can be solved.

In the present article we
solve differential equation (\ref{cp}) for   a large class of homogeneous spaces $G/H$ under a certain algebraic condition on the tangent space $T_{eH}(G/H)$.  Firstly, we recall the following general context.

Let $M=G/H$ be a homogeneous space and let $\fr{g},\fr{h}$ be the Lie algebras of $G,H$ respectively.  Consider a reductive decomposition $\fr{g}=\fr{h}\oplus \fr{m}$, $\operatorname{Ad}(H)\fr{m}\subseteq \fr{m}$, with respect to an $\operatorname{Ad}$-invariant inner product $B$ of $\fr{g}$, where $\operatorname{Ad}:G\rightarrow \operatorname{Aut}(\fr{g})$ denotes the adjoint representation of $G$.  Such a product exists if and only if $G$ is the direct product of a compact group with a vector group diffeomorphic to $\mathbb R^n$ (\cite{Mi}).  An example of such a product is the negative of the Killing form in simple Lie algebras.  Moreover, assume that the Lie algebra $\fr{h}$ has non trivial center $\fr{z}(\fr{h})$.  The tangent space $T_oM$ at $o=eH$ can be identified to $\fr{m}$.
Let $\pi: G\rightarrow G/H$ be the projection and for $p\in G$, let $\tau_p:G/H\rightarrow G/H$ be the left translation in $G/H$ by $p$.  Any $G$-invariant Riemannian metric in $M$, i.e. any metric invariant by the left translations, corresponds to a unique $\operatorname{Ad}(H)$-invariant inner product $\langle \ ,\ \rangle$ in $\fr{m}$, and in turn to a unique linear endomorphism $A:\fr{m}\rightarrow \fr{m}$ satisfying 

\begin{equation}\label{metend}\langle X,Y\rangle=B(AX,Y),\quad X,Y\in \fr{m}.\end{equation}  

\noindent The endomorphism $A$ is symmetric with respect to $B$, positive definite and $\operatorname{Ad}(H)$-equivariant, that is $\operatorname{Ad}(h)\circ A=A\circ \operatorname{Ad}(h)$, $h\in H$.  Hence, the endomorphism $A$ is diagonalizable and its eigenspaces are $\operatorname{Ad}(H)$-invariant and $B$-orthogonal.

Choose $W\in \fr{z}(\fr{h})$ and let $\omega^W$ be the $\operatorname{Ad}(H)$-invariant form on $G/H$, which on the origin $o$ of $G/H$ is defined by 

\begin{equation}\label{omega}\omega^W_o(X,Y)=-B(W,[X,Y]),\quad X,Y\in \fr{m}.\end{equation}

 \noindent For example, when $G/H$ is an adjoint orbit, $\omega^W$ is induced from the Kirillov-Konstant-Souriau form (see \cite{Bo-Jo}).  Moreover, let $I$ be the $G$-invariant $(1,1)$ tensor such that

\begin{equation}\label{tensor}\omega^W(X,Y)=g(X,IY).\end{equation}

\noindent  Then $I$ satisfies the equation $g(IX,Y)+g(Y,IX)=0$ for all vector fields $X,Y$.\\

We assume that the following condition is satisfied:\\

\noindent \textbf{Condition A.} There exist eigenspaces $\fr{m}_a,\fr{m}_b\subseteq \fr{m}$ of the metric endomorphism $A$, with corresponding eigenvalues $\lambda_a,\lambda_b$, such that

\begin{equation}\label{mod}[\fr{m}_a,\fr{m}_b]\subseteq \fr{m}_a.\end{equation}

We have the following result.

\begin{theorem}\label{maintheorem}
Let $(G/H, g)$ be a Riemannian homogeneous space, with $G$ compact and $H$ having non-discrete center, satisfying Condition A.  Let $x(t)$ be the motion of a charged particle under electromagnetic field $kI$, namely $\nabla_{\dot{x}}{\dot{x}}=kI(\dot{x})$, where $I$ is given by Equation (\ref{tensor}).  If $x$ satisfies the initial conditions

\begin{equation}\label{in}x(0)=o \quad \makebox{and} \quad \dot{x}(0)\in \fr{m}_a\oplus \fr{m}_b,
\end{equation}

\noindent then it is given by

\begin{equation*}\label{curve}x(t)=\exp t(X_a+\frac{\lambda_b}{\lambda_a} X_b+\frac{k}{\lambda_a}W) \exp t(1-\frac{\lambda_b}{\lambda_a})(X_b+\frac{k}{\lambda_b}W)\cdot o,
\end{equation*}

\noindent where $X_a,X_b$ are the projections of $\dot{x}(0)$ on $\fr{m}_a,\fr{m}_b$ respectively. 
\end{theorem}

As a consequence, we obtain the following description of the corresponding charged particle motion in homogeneous fibrations $p:G/H\rightarrow G/K$, where the normal metric on $G/H$ is deformed along the fiber $K/H$: 

\begin{corol}\label{corol}Let $G/H$ be a homogeneous space with $G$ compact and $H$ having non-discrete center.  Consider a homogeneous fibration $p:G/H\rightarrow G/K$ with fiber $K/H$, where $K$ is a closed Lie subgroup of $G$, so that $T_{o}(G/H)=T_{p(o)}(G/K)\oplus T_o(K/H)$.  Endow $G/H$ with the $G$-invariant metric $g$ induced from the inner product 

\begin{equation}\label{Ch}\langle \ ,\ \rangle=\left.B\right|_{T_{p(o)}(G/K)\times T_{p(o)}(G/K)}+\lambda\left.B\right|_{T_o(K/H)\times T_o(K/H)},\end{equation}

\noindent where $B$ is an $\operatorname{Ad}$-invariant inner product in $\fr{g}$.  Let $x(t)$ be the motion of a charged particle through the origin $o$, under electromagnetic field $kI$ where $I$ is given by Equation (\ref{tensor}).  Then $x(t)$ is given by the equation   

\begin{equation}\label{curve1}x(t)=\exp t(X_a+\lambda X_b+kW) \exp t(1-\lambda)(X_b+\frac{k}{\lambda}W)\cdot o,
\end{equation}

\noindent where $X_a,X_b$ are the projections of $\dot{x}(0)$ on $\fr{m}_a=T_{p(o)}(G/K)$ and $\fr{m}_b=T_o(K/H)$ respectively. 
\end{corol}

Metrics of the form (\ref{Ch}) are examples of Cheeger deformation metrics (\cite{Che}) and are widely used in geometry.  Moreover, particular examples of homogeneous fibrations, such as the Hopf fibration of $\mathbb S^3$, are implemented in gauge theory and its applications to particle physics (\cite{BMSS}).\\
  \indent The above Theorem generalizes Ikawa's main result (\cite[Theorem 1.1]{Ik3}) in the sense that it describes the motion of a charged particle for all metrics induced by a homogeneous fibration of a compact homogeneous space.

\begin{example}(Twistor fibration on generalised flag manifolds). 
For $G$ semisimple let $M=G/H$ be a generalised flag manifold, i.e. an adjoint orbit of an element $W$ in  $\fr{g}$. It is known  (\cite{Bu-Ra}) that any flag manifold $M=G/H$ can 
be fibered over a compact inner symmetric space $G/K$ ($H\subset K$) under the twistor fibration
$p:G/H\rightarrow G/K$.
The {\it normal metric} of $G/H$ 
is the $G$-invariant metric induced by the negative of the Killing form  $\fr{g}$, denoted by $B$.  
We endow $G/H$ with the ``deformation" of the normal metric
along the fibers $K/H$ of the twistor fibration , given by

\begin{equation*}\langle \ ,\ \rangle=\left.B\right|_{T_{p(o)}(G/K)\times T_{p(o)}(G/K)}+\lambda\left.B\right|_{T_o(K/H)\times T_o(K/H)},\end{equation*}
 
\noindent 

Then, by virtue of Corollary \ref{corol}, the equation of motion of a charged particle in $(M,\langle \ ,\ \rangle)$ under electromagnetic field $kI$, is given by (\ref{curve1}).

\end{example}

\section{Proof of the main results}

\subsection{Proof of Theorem \ref{maintheorem}}

 To simplify the exposition and since $G$ is compact, we may assume that $G$ is a matrix group so that the left and right translations $L_p,R_p$ in $G$ are matrix multiplications their respective differentials are also given by matrix multiplications.  In particular, for $p,q\in G$ and $v\in T_qG$ we have $L_p(q)=pq$, $(dL_p)_qv=pv$ and similarly for $R_q$.  Moreover, for $\pi(q)\in G/H$ we denote by $p\pi(q)$ the left translation $\tau_p(\pi(q))$ in $G/H$, while for $w\in T_{\pi(q)}(G/H)$ we denote by $pw$ the differential $(d\tau_p)_{\pi(q)}w$.  Under the above notations, we also recall the identities $p\pi(q)=\pi(pq)$ and $p(d\pi_qv)=d\pi_{pq}(pv)$.\\  
\indent For $Z\in \fr{g}$, we set $e^Z:=\exp(Z)$.  We also set 

\begin{equation}\label{relxy}X:=X_a+\frac{\lambda_b}{\lambda_a} X_b+ \frac{k}{\lambda_a}W,\quad Y:=(1-\frac{\lambda_b}{\lambda_a})(X_b+\frac{k}{\lambda_b}W),\quad \makebox{and} \quad \alpha(t):=e^{tX}e^{tY}.\end{equation} 

\noindent  We will show that $x(t):=\pi(\alpha(t))$ satisfies Equation (\ref{cp}).  Equivalently, we will show that $g(V,\nabla_{\dot{x}}{\dot{x}})=g(V,kI(\dot{x}))$ for any vector field $V$ in $G/H$.  However, we may substitute $V$ with any Killing vector field $Z^*$, $Z\in \fr{g}$, given by

\begin{equation}\label{kil}Z^*_{\pi(p)}=\left.\frac{d}{dt}\right|_{t=0}\pi(\exp(tZ)p)=d\pi_{p}(Zp),\quad \pi(p)\in G/H,\quad p\in G,\end{equation}

\noindent as the above vector fields generate the tangent bundle $T(G/H)$.  Hence, we need to show that
\begin{equation}\label{have}g(Z^*, \nabla_{\dot{x}}{\dot{x}})=g(Z^*,kI(\dot{x})),\end{equation}

\noindent for any vector $Z$ in $\fr{g}$.\\
\indent  The derivative of $x(t)=\pi(e^{tX}e^{tY})=\pi(\alpha(t))$ is $\dot{x}(t)=d\pi_{\alpha(t)}(\dot{\alpha}(t))$.  We will express $\dot{\alpha}(t)$ explicitly. Given that $\left.\frac{d}{dt}e^{tZ}\right.=Ze^{tZ}=e^{tZ}Z$ and $\operatorname{Ad}(e^{tZ})Z=e^{tZ}Ze^{-tZ}=Z$ for all $Z\in \fr{g}$ and $t\in \mathbb R$, we have

\begin{eqnarray*}\dot{\alpha}(t)&=&\left.\frac{d}{dt}(e^{tX}e^{tY})\right.=(\left.\frac{d}{dt}e^{tX}\right.)e^{tY}+e^{tX}(\left.\frac{d}{dt}e^{tY}\right.)=Xe^{tX}e^{tY}+e^{tX}e^{tY}Y=X\alpha(t)+\alpha(t)Y\nonumber \\
&=&
\alpha(t)\big(\alpha^{-1}(t)X\alpha(t)+Y\big)=\alpha(t)\big(e^{-tY}e^{-tX}Xe^{tX}e^{tY}+Y\big)=\alpha(t)\big(e^{-tY}Xe^{tY}+e^{-tY}Ye^{tY}\big)\nonumber \\
&=&
\alpha(t)\big(\operatorname{Ad}(e^{-tY})(X+Y)\big).
\end{eqnarray*}
    
\noindent Using the above expression for $\dot{\alpha}(t)$, we obtain

\begin{equation}\label{derx}\dot{x}(t)=d\pi_{\alpha(t)}(\dot{\alpha}(t))=d\pi_{\alpha(t)}\bigg(\alpha(t)\big(\operatorname{Ad}(e^{-tY})(X+Y)\big)\bigg).\end{equation}

\indent Next, we will express explicitly the left-hand side of Equation (\ref{have}).  By using Koszul's formula (\cite{On} p. 61, we have 

\begin{equation} \label{koz1} g(Z^*,\nabla_{\dot{x}}{\dot{x}})=\dot{x}g(Z^*,\dot{x})+g(\dot{x},[Z^*,\dot{x}])-\frac{1}{2}Z^*g(\dot{x},\dot{x}).\end{equation}

By using the compatibility of $\nabla$ with the metric and the fact that it is torsion-free (see for example \cite{On}, p. 61), we obtain $Z^*g( \dot{x},\dot{x})=2g( \nabla_{Z^*}{\dot{x}},\dot{x})$ and $\nabla_{Z^*}{\dot{x}}-\nabla_{\dot{x}}{Z^*}=[Z^*,\dot{x}]$.  On the other hand, as $Z^*$ is Killing, we have $( \nabla_{\dot{\gamma}}{Z^*},\dot{\gamma})=0$ (\cite{On}, p. 251).  The last three equations imply that $g(\dot{x},[Z^*,\dot{x}])-\frac{1}{2}Z^*g(\dot{x},\dot{x})=0$, hence by virtue of Equation (\ref{koz1}) we have

\begin{equation*}\label{koz2}g(Z^*,\nabla_{\dot{x}}{\dot{x}})=\dot{x}g(Z^*,\dot{x}).\end{equation*}

By taking into account the invariance of $g$ by the left-translations, by substituting Equations (\ref{derx}) and (\ref{kil}) into the above equation, using the facts that $x(t)=\pi(\alpha(t))$ and $\alpha^{-1}(t)\dot{\alpha}(t)=\operatorname{Ad}(e^{-tY})(X+Y)$ as well as the identification $\fr{m}=d\pi_e(\fr{g})$, we obtain the following

\begin{eqnarray}g(Z^*,\nabla_{\dot{x}}{\dot{x}})(t)&=&\dot{x}g(Z^*,\dot{x})(t)=\frac{d}{dt}g(Z^*_{x(t)},\dot{x}(t))=\frac{d}{dt}g\big(d\pi_{\alpha(t)}(Z\alpha(t)),d\pi_{\alpha(t)}(\dot{\alpha}(t))\big)\nonumber\\
&=&
\frac{d}{dt}g_o\big(\alpha^{-1}(t)d\pi_{\alpha(t)}(Z\alpha(t)),\alpha^{-1}(t)d\pi_{\alpha(t)}(\dot{\alpha}(t))\big)\nonumber \\
&=&
\frac{d}{dt}\big\langle d\pi_{e}(\alpha^{-1}(t)Z\alpha(t)),d\pi_{e}(\alpha^{-1}(t)\dot{\alpha}(t))\big\rangle\nonumber \\
&=&
\frac{d}{dt}\big\langle \big(\operatorname{Ad}(\alpha^{-1}(t))Z\big)_{\fr{m}},\big(\operatorname{Ad}(e^{-tY})(X+Y)\big)_{\fr{m}}\big\rangle\label{lh1}
\end{eqnarray}

\noindent (here the subscript $\fr{m}$ denotes the projection on $\fr{m}$).  We abuse the notation $\dot{\alpha}^{-1}(t)$ to denote the derivative $\left.\frac{d}{dt}\alpha^{-1}(t)\right.$, and we differentiate the relation $\alpha^{-1}(t)\alpha(t)=e$ to obtain $\dot{\alpha}^{-1}(t)\alpha(t)=-\alpha^{-1}(t)\dot{\alpha}(t)=-\operatorname{Ad}(e^{-tY})(X+Y)$. Also, we set

\begin{equation}\label{f}f(t):=\operatorname{Ad}(e^{-tY})(X+Y)=\alpha^{-1}(t)\dot{\alpha}(t) \quad \makebox{and} \quad \widetilde{Z}(t):=\operatorname{Ad}(\alpha^{-1}(t))Z.\end{equation}

\noindent Then 

\begin{equation*}\left.\frac{d}{dt}f(t)\right.=[f(t),Y] \ \ \makebox{and}\end{equation*}

\begin{eqnarray*}\left.\frac{d}{dt}\widetilde{Z}(t)\right.&=& \dot{\alpha}^{-1}(t)Z\alpha(t)+\alpha^{-1}(t)Z\dot{\alpha}(t)=\dot{\alpha}^{-1}(t)\alpha(t)\widetilde{Z}(t)+\widetilde{Z}(t)\alpha^{-1}(t)\dot{\alpha}(t)\\
&=&
-\alpha^{-1}(t)\dot{\alpha}(t)\widetilde{Z}(t)+\widetilde{Z}(t)\alpha^{-1}(t)\dot{\alpha}(t)=[\widetilde{Z}(t),f(t)].
\end{eqnarray*}

Then Equation (\ref{lh1}) yields

\begin{equation}\label{lh2}g(Z^*,\nabla_{\dot{x}}{\dot{x}})(t)=\frac{d}{dt}\big\langle \widetilde{Z}(t)_{\fr{m}},f(t)_{\fr{m}} \big\rangle=\big\langle[\widetilde{Z}(t),f(t)]_{\fr{m}},f(t)_{\fr{m}}\big\rangle+\big\langle\widetilde{Z}(t)_{\fr{m}},[f(t),Y]_{\fr{m}}\big\rangle \end{equation}

\noindent Along with Equation (\ref{metend}), the mutual $B$-orthogonality of $\fr{m}$ and $\fr{h}$ as well as the $\operatorname{Ad}$-invariance of $B$, Equation (\ref{lh2}) yields

\begin{eqnarray}g(Z^*,\nabla_{\dot{x}}{\dot{x}})(t) &=& B\big([\widetilde{Z}(t),f(t)]_{\fr{m}},Af(t)_{\fr{m}}\big)+B\big(\widetilde{Z}(t)_{\fr{m}},A[f(t),Y]_{\fr{m}}\big)\nonumber \\
&=&
B\big([\widetilde{Z}(t),f(t)],Af(t)_{\fr{m}}\big)+B\big(\widetilde{Z}(t),A[f(t),Y]_{\fr{m}}\big)\nonumber \\
&=&
B\big(\widetilde{Z}(t),[f(t),Af(t)_{\fr{m}}]+A[f(t),Y]_{\fr{m}}\big).\label{13}
\label{lh3}\end{eqnarray}

 We will express the term $[f(t),Af(t)_{\fr{m}}]+A[f(t),Y]_{\fr{m}}$ explicitly.  Firstly, recall that $\operatorname{Ad}(e^v)w=e^{\operatorname{ad}(v)}w$ for all $v,w\in \fr{g}$.  Using this fact, along with relation (\ref{mod}), we deduce that

\begin{equation}\label{incl}\operatorname{Ad}(e^{\fr{m}_b})\fr{m}_a\subseteq \fr{m}_a.\end{equation}

\noindent By relation (\ref{relxy}), we have $Y:=(1-\frac{\lambda_b}{\lambda_a})(X_b+\frac{k}{\lambda_b}W)$, and thus $\operatorname{Ad}(e^{-tY})(X_b+\frac{k}{\lambda_b}W)=(X_b+\frac{k}{\lambda_b}W)$.  Moreover, relation (\ref{incl}) and the $\operatorname{ad}_{\fr{h}}$-invariance of $\fr{m}_a$ imply that $\operatorname{Ad}(e^{\fr{m}_b\oplus \fr{h}})\fr{m}_a\subseteq \fr{m}_a$, and thus $\operatorname{Ad}(e^{-tY})(X_a)\in \fr{m}_a$.  Using the above facts, along with relations (\ref{relxy}) and (\ref{f}), we obtain 

\begin{eqnarray}f(t)&=&\operatorname{Ad}(e^{-tY})(X+Y)=\operatorname{Ad}(e^{-tY})(X_a+X_b+\frac{k}{\lambda_b}W)\nonumber \\
&=&
\operatorname{Ad}(e^{-tY})X_a+\operatorname{Ad}(e^{-tY})(X_b+\frac{k}{\lambda_b}W)\nonumber\\
&=&
\underbrace{\operatorname{Ad}(e^{-tY})X_a}_{f(t)_{\fr{m}_a}}+\underbrace{X_b}_{f(t)_{\fr{m}_b}}+\underbrace{\frac{k}{\lambda_b}W}_{f(t)_{\fr{h}}}\label{proj}
\end{eqnarray}

\noindent (here the subscripts denote the projections of $f(t)$ on the respective subspaces). Taking into account relation (\ref{mod}) and the $\operatorname{ad}(\fr{h})$-invariance of the spaces $\fr{m}_a,\fr{m}_b$, we obtain $[f(t)_{\fr{m}_a},f(t)_{\fr{m}_b}]\in \fr{m}_a$ and $[f(t)_{\fr{m}_a},f(t)_{\fr{h}}]\in \fr{m}_a$. Along with Equation (\ref{proj}), the fact that $Y=(1-\frac{\lambda_b}{\lambda_a})(f_{\fr{m}_b}+f_{\fr{h}})$, the fact that $\fr{m}_a,\fr{m}_b$ are eigenspaces of $A:\fr{m}\rightarrow \fr{m}$ with corresponding eigenvalues $\lambda_a,\lambda_b$, and omitting $t$ for simplicity, we deduce that

\begin{eqnarray*}[f,Af_{\fr{m}}]+A[f,Y]_{\fr{m}}&=&[f,Af_{\fr{m}_a}+Af_{\fr{m}_b}]+(1-\frac{\lambda_b}{\lambda_a})A[f,f_{\fr{m}_b}+f_{\fr{h}}]_{\fr{m}}\nonumber \\
&=&
[f,Af_{\fr{m}_a}+Af_{\fr{m}_b}]+(1-\frac{\lambda_b}{\lambda_a})\big(A[f_{\fr{m}_a},f_{\fr{m}_b}]_{\fr{m}}+A[f_{\fr{m}_a},f_{\fr{h}}]_{\fr{m}}\big) \nonumber \\
&=&
[f,\lambda_a f_{\fr{m}_a}+\lambda_b f_{\fr{m}_b}]+(1-\frac{\lambda_b}{\lambda_a})\big(\lambda_a[f_{\fr{m}_a},f_{\fr{m}_b}]+\lambda_a[f_{\fr{m}_a},f_{\fr{h}}]
\big)\nonumber \\
&=&
(\lambda_b-\lambda_a)[f_{\fr{m}_a},f_{\fr{m}_b}]+\lambda_a [f_{\fr{h}},f_{\fr{m}_a}]+\lambda_b [f_{\fr{h}},f_{\fr{m}_b}]+(\lambda_a-\lambda_b)\big([f_{\fr{m}_a},f_{\fr{m}_b}]+[f_{\fr{m}_a},f_{\fr{h}}]\big)\nonumber \\
&=&
\lambda_b\big([f_{\fr{h}},f_{\fr{m}_a}]+[f_{\fr{h}},f_{\fr{m}_b}]\big)=\lambda_b[f_{\fr{h}},f_{\fr{m}}]=k[W,f_{\fr{m}}].
\end{eqnarray*}

Therefore, by virtue of Equation (\ref{13}) we deduce that the left-hand side of Equation (\ref{have}) becomes 

\begin{equation}\label{LH}g(Z^*, \nabla_{\dot{x}}{\dot{x}})(t)=kB\big(\widetilde{Z}(t),[W,f(t)_{\fr{m}}]\big).\end{equation}

To conclude the proof it remains to compute the right-hand side of Equation (\ref{have}) and compare it with Equation (\ref{LH}).  Firstly, we find the explicit form of the tensor $I$.  Since it is left-invariant, it suffices to determine $I_o:\fr{m}\rightarrow \fr{m}$.  From Equations (\ref{tensor}), (\ref{omega}) and (\ref{metend}) we obtain for any $X,Y\in \fr{m}$ that $\langle I_oX,Y\rangle =-\langle X,I_oY\rangle =-g_o(X,I_oY)=-\omega_o^W(X,Y)=B(W,[X,Y])=B([W,X],Y)=\langle A^{-1}[W,X],Y\rangle$, therefore 

\begin{equation*}I_oX=A^{-1}[W,X].\end{equation*}

Using the above equation, the $G$-invariance of the metric and the tensor $I$, relation (\ref{kil}), Equation (\ref{derx}), relations (\ref{f}), Equation (\ref{metend}), the $\operatorname{ad}(\fr{h})$-invariance of $\fr{m}$ and the $B$-orthogonality of $\fr{m}$ and $\fr{h}$, we deduce that the right-hand side of Equation (\ref{have}) becomes

\begin{eqnarray*}g(Z^*,kI(\dot{x}))(t)&=&k\big\langle \alpha^{-1}(t)d\pi_{\alpha(t)}(Z\alpha(t)),I_o(\alpha^{-1}(t)\dot{x}(t))\big\rangle\\
&=&
k\big\langle d\pi_{e}(\alpha^{-1}(t)Z\alpha(t)),I_o\big(\alpha^{-1}(t)d\pi_{\alpha(t)}(\dot{\alpha}(t))\big)\big\rangle\\
&=&
k\big\langle d\pi_{e}(\operatorname{Ad}(\alpha^{-1}(t))Z),I_o\big(d\pi_{e}(\alpha^{-1}(t)\dot{\alpha}(t))\big)\big\rangle\\
&=&
k\big\langle \big(\operatorname{Ad}(\alpha^{-1}(t))Z\big)_{\fr{m}},I_o(f(t)_{\fr{m}})\big\rangle=k\big\langle \widetilde{Z}(t)_{\fr{m}},A^{-1}[W,f(t)_{\fr{m}}]\big\rangle\\
&=&
kB\big( \widetilde{Z}(t)_{\fr{m}},[W,f(t)_{\fr{m}}]\big)=kB\big( \widetilde{Z}(t),[W,f(t)_{\fr{m}}]\big),
\end{eqnarray*}
  
\noindent which by virtue of Equation (\ref{LH}), concludes the proof.\qed

\subsection{Proof of Corollary \ref{corol}}

Let $\fr{g},\fr{k},\fr{h}$ be the Lie algebras of the groups $G,K,H$ respectively. The spaces $\fr{m}_a=T_{p(o)}(G/K)$ and $\fr{m}_b=T_o(K/H)$ arise from the reductive decompositions $\fr{g}=\fr{k}\oplus \fr{m}_a$ and $\fr{k}=\fr{h}\oplus \fr{m}_b$.  The first reductive decomposition implies that $[\fr{k},\fr{m}_a]\subseteq \fr{m}_a$ and thus $[\fr{m}_a,\fr{m}_b]\subseteq \fr{m}_a$.  Corollary \ref{corol} then follows from Theorem \ref{maintheorem}, with $\lambda_a=1$, $\lambda_b=\lambda$, along with the fact that $T_o(G/H)=\fr{m}_a\oplus \fr{m}_b$.\qed

\end{document}